\documentclass[a4paper,12pt]{article}
\usepackage{amssymb,amsmath,bm}
\usepackage{graphicx}
\usepackage{enumerate}

\title{Differential Characterization of Quasi-Concave Functions without Twice Differentiability}
\author{Yuhki Hosoya\thanks{TEL: +81-90-5525-5142, E-mail: ukki(at)gs.econ.keio.ac.jp}\\ Faculty of Economics, Chuo University\thanks{742-1, Higashinakano, Hachioji-shi, Tokyo, 192-0393, Japan.}}
\date{}
\begin{document}
\maketitle

\begin{abstract}
This paper presents a necessary and sufficient condition for a real-valued function defined on an open and convex subset of a Banach space to be quasi-concave, and a sufficient condition for such a function to be strictly quasi-concave. These conditions are applicable to continuously differentiable functions that satisfy a mild additional assumption, and do not require the functions to be twice differentiable. Because this additional assumption is trivially satisfied for twice continuously differentiable functions, our results are pure extensions to classical results.

\vspace{12pt}

\vspace{12pt}

\vspace{12pt}
\end{abstract}

\section{Introduction}
The notion of quasi-concave functions is frequently used for optimization problems in economics. However, it is difficult to characterize the quasi-concavity of a function in terms of conditions for its derivatives. For concavity, in contrast, there is a clear characterization: for example, a continuously differentiable real-valued function defined on an open and convex subset of the real line is concave if and only if its derivative is nonincreasing. However, for quasi-concavity, such a clear characterization is absent.

Otani (1983) treated this problem and provided a necessary and sufficient condition for twice continuously differentiable non-degenerate functions defined on an open and convex subset in $\mathbb{R}^n$ to be quasi-concave. His condition is that the Hessian matrix of this function is negative semi-definite on the kernel of the gradient vector. However, there are several optimization problems in which twice continuous differentiability of the objective function cannot be assumed. The main purpose of this paper is to extend Otani's result to functions that are not necessarily twice continuously differentiable but are once continuously differentiable and satisfy an additional property. Note that, any twice continuously differentiable function satisfies this additional property, and so our result is a pure extension to a known result.

Moreover, we extend some classical results for quasi-concave functions on $\mathbb{R}^n$ to any Banach space. Specifically, we treat a continuously differentiable real-valued function defined on an open and convex subset of a Banach space, and present three results. First, we provide a necessary and sufficient condition for a non-degenerate function to be quasi-concave (Theorem 1). Second, we provide a sufficient condition for a function to be strictly quasi-concave (Theorem 2). As we noted above, both results are pure extensions to classical known results.

Note that our ``additional property'' introduced in this paper is deeply related to Frobenius' theorem. Debreu (1972) considered the following {\bf total differential equation}:
\[\nabla f(x)=\lambda(x)g(x),\]
where $g(x)$ is the given function and the pair of $f(x)$ and $\lambda(x)$ is a solution. In Debreu's theory, the vector field $g(x)$ exhibits the price vector under which consumption plan $x$ is chosen. Thus, the above equation concerns Lagrange's multiplier rule, and $f(x)$ can be treated as the objective function for the consumer. Debreu expected that if $g(x)$ is continuously differentiable, then the solution $f(x)$ must be twice continuously differentiable. However, a counterexample found by two students was later presented by Debreu (1976), in which the solution $f(x)$ must not be twice continuously differentiable. Therefore, we cannot assume that the objective function for the consumer is twice continuously differentiable if $g(x)$ is only continuously differentiable. Our additional property implies the continuous differentiability for such a $g(x)$, and so our results are applicable to consumer theory.

In Subsection 2.1, we introduce some knowledge on Banach spaces that is needed to understand our results. In Subsection 2.2, we present the main results of this study. Section 3 contains several remarks concerning our results.

\section{Main Results}

\subsection{Preliminaries}
In this subsection, we introduce some basic knowledge on Banach spaces. We assume that readers know the definition of a Banach space. All of the facts introduced in this subsection are explained and proved in standard textbooks: see, for example, chapters 2-5 of Luenberger (1969) or chapters 5-6 of Aliprantis and Border (2006).

Suppose that $X, Y$ are Banach spaces. Let $L(X,Y)$ be the set of all linear and continuous functions from $X$ into $Y$. Then, $L(X,Y)$ becomes a Banach space with respect to the following {\bf operator norm}
\[\|A\|=\sup\{\|Ax\||x\in X,\ \|x\|\le 1\}.\]
In particular, the space $X'$ of all linear and continuous functions from $X$ into $\mathbb{R}$ is called the {\bf dual space} of $X$. Note that $X'$ is the same as $L(X,\mathbb{R})$, and thus it is also a Banach space with the operator norm. For $x'\in X'$, we write $\langle x,x'\rangle$ instead of $x'(x)$. It is easy to show that the function $f(x,x')=\langle x,x'\rangle$ is continuous.

If $X=\mathbb{R}^n$, then $X'=\mathbb{R}^n$ and $\langle x,x'\rangle$ is the usual inner product between two vectors $x=(x_1,...,x_n)$ and $x'=(x_1',...,x_n')$; that is,
\[\langle x,x'\rangle=\sum_{i=1}^nx_ix_i'.\]
In this paper, $A^T$ denotes the transpose of an $n\times m$ matrix $A$. Note that, if $X=\mathbb{R}^n$, then $\langle x,x'\rangle$ can also be written as $x^Tx'$.

Next, suppose that $X,Y$ are Banach spaces. Suppose that $U$ is an open set in $X$ and $f:U\to Y$. If there exists an element $A\in L(X,Y)$ such that
\[\lim_{h\to 0}\frac{\|f(x+h)-f(x)-Ah\|}{\|h\|}=0,\]
then $f$ is said to be {\bf Fr\'echet differentiable} at $x$, and $A$ is called the {\bf Fr\'echet derivative} of $f$ at $x$. It is known that if $f$ is Fr\'echet differentiable at $x$, then it is continuous at $x$, and its Fr\'echet derivative at $x$ is unique. In this paper, $Df(x)$ denotes the Fr\'echet derivative of $f$ at $x$. If $f$ is Fr\'echet differentiable at any point of $U$ and $Df:x\mapsto Df(x)$ is continuous, then we say that $f$ is $C^1$. We often abbreviate ``Fr\'echet differentiable'' to ``differentiable'' simply.

Specifically, suppose that $Y=\mathbb{R}$. Let $U$ be an open subset of the Banach space $X$ and $f:U\to \mathbb{R}$ be $C^1$. Then, the mapping $Df:x\mapsto Df(x)$ is a function from $U$ into $X'$, and $X'$ is a Banach space. Therefore, we can define the Fr\'echet derivative of $Df$ at $x$, which is denoted by $D^2f(x)$. If $D^2f(x)$ is defined on $U$ and $D^2f:x\mapsto D^2f(x)$ is continuous, then we say that $f$ is $C^2$.

Suppose that $X, Y, Z$ are Banach spaces, $U\subset X$ and $V\subset Y$ are open sets, $f:V\to Z$ and $g:U\to Y$ are given, $g(x)\in V$, $g$ is differentiable at $x$, and $f$ is differentiable at $g(x)$. Then, $f\circ g$ is differentiable at $x$, and
\[D(f\circ g)(x)=Df(g(x))\circ Dg(x).\]
This formula is called the {\bf chain rule}.

Suppose that $X$ is a linear space and $U$ is a convex subset of $X$. A function $f:U\to \mathbb{R}$ is said to be {\bf quasi-concave} if and only if for every $x,y\in U$ and $t\in [0,1]$,
\[f((1-t)x+ty)\ge \min\{f(x),f(y)\}.\]
If, in addition,
\[f((1-t)x+ty)>\min\{f(x),f(y)\}\]
when $x\neq y$ and $0<t<1$, then we say that $f$ is {\bf strictly quasi-concave}.

\subsection{Results}
Suppose that $X$ is a Banach space and $U$ is an open subset of $X$. For $f:U\to \mathbb{R}$, we say that $f$ is $C^1_*$ if $f$ is $C^1$ and there exists a pair of a $C^1$ function $g:U\to X'$ and a positive continuous function $\lambda:U\to \mathbb{R}$ such that
\begin{equation}\label{VEC}
Df(x)=\lambda(x)g(x).
\end{equation}
Clearly, if $f$ is $C^2$, then $f$ is $C^1_*$: choose $\lambda\equiv 1$ and $g(x)=Df(x)$. When $X$ is a Hilbert space, then the norm of $X$ is differentiable everywhere except the zero vector. In this case, if $f$ is $C^1$ and $Df(x)\neq 0$ for all $x\in U$, then $f$ is $C^1_*$ if and only if the following function
\[h(x)=\frac{1}{\|Df(x)\|}Df(x)\]
is $C^1$. Indeed, if $h$ is $C^1$, then we can choose $\lambda(x)=\|Df(x)\|$ and $g(x)=h(x)$. Conversely, if there exists a $C^1$ function $g$ and a continuous function $\lambda$ such that (\ref{VEC}) holds, then $g(x)\neq 0$ for all $x\in U$, and
\[h(x)=\frac{1}{\|g(x)\|}g(x),\]
which is $C^1$. In particular, if $X=\mathbb{R}^n$, then the space of all $C^1_*$ functions whose derivatives never vanish coincides with the space of $C^1$ functions $f$ such that the above $h$ is $C^1$.

Our first main result is as follows.

\vspace{12pt}
\noindent
{\bf Theorem 1}. Suppose that $U$ is an open and convex subset of a Banach space $X$, $f:U\to \mathbb{R}$ is $C^1_*$, and $Df(x)\neq 0$ for all $x\in U$. Choose a pair of a $C^1$ function $g$ and a positive continuous function $\lambda$ that satisfies (\ref{VEC}). Then, the following two claims are equivalent.

\begin{enumerate}[1)]
\item $f$ is quasi-concave.

\item $\langle w, Dg(x)w\rangle \le 0$ for any $x\in U$ and $w\in X$ such that $\langle w,g(x)\rangle=0$.
\end{enumerate}

\vspace{12pt}
\noindent
{\bf Proof}. First, we prove a lemma.

\vspace{12pt}
\noindent
{\bf Lemma 1}(Lagrange's multiplier rule). Suppose that $U$ is an open and convex subset of a Banach space $X$ and $f:U\to \mathbb{R}$ is continuous and quasi-concave. Choose any $x'\in X'\setminus \{0\}$. If $f$ is differentiable at $x^*$ and $Df(x^*)=\lambda x'$ for some $\lambda>0$, then $f(x^*)\ge f(x)$ for every $x\in U$ such that $\langle x,x'\rangle=\langle x^*,x'\rangle$.

\vspace{12pt}
\noindent
{\bf Proof of Lemma 1}. Suppose not. Then, there exists $x\in U$ such that $\langle x,x'\rangle=\langle x^*,x'\rangle$ and $f(x)>f(x^*)$. Because $f$ is continuous and $U$ is open, there exists $y\in U$ such that $\langle y,x'\rangle<\langle x^*,x'\rangle$ and $f(y)>f(x^*)$. Define $y(t)=(1-t)x^*+ty$. Then, $y(t)\in U$ for all $t\in [0,1]$, and by the quasi-concavity of $f$, we have that $f(y(t))\ge f(x^*)$ for all $t\in [0,1]$. Therefore,
\[0>\langle y-x^*,Df(x^*)\rangle=\left.\frac{d}{dt}f(y(t))\right|_{t=0}=\lim_{t\downarrow 0}\frac{f(y(t))-f(x^*)}{t}\ge 0,\]
which is a contradiction. This completes the proof. $\blacksquare$

\vspace{12pt}
Suppose that 1) holds. Choose any $x\in U$ and $w\in X$ such that $\langle w,g(x)\rangle=0$. Define
\[c(t)=f(x+tw).\]
Then, there exists $\varepsilon>0$ such that $c(t)$ is defined on $I=[-\varepsilon,\varepsilon]$. By Lagrange's multiplier rule, we have that $c(t)$ attains the maximum at $t=0$. Hence, $c(t)\le c(0)$ for all $t\in [0,\varepsilon]$. By the mean value theorem, we can obtain a sequence $(t_k)$ such that $t_k\downarrow 0$ as $k\to \infty$ and $c'(t_k)\le 0$ for all $k$. Therefore,
\begin{align*}
0\ge&~\limsup_{k\to \infty}\frac{c'(t_k)}{t_k}=\limsup_{k\to \infty}\frac{\langle w,Df(x+t_kw)\rangle}{t_k}\\
=&~\limsup_{k\to \infty}\frac{\lambda(x+t_kw)\langle w,g(x+t_kw)\rangle}{t_k}\\
\ge&~M\limsup_{k\to \infty}\frac{\langle w,g(x+t_kw)\rangle}{t_k}\\
=&~M\limsup_{k\to \infty}\frac{\langle w,g(x+t_kw)-g(x)\rangle}{t_k}=M\langle w,Dg(x)w\rangle,
\end{align*}
where $M=\max\{\lambda(x+tw)|t\in I\}>0$, which implies that 2) holds. Therefore, 1) implies 2).

Next, suppose that 2) holds and 1) is violated. Then, there exist $x,y\in U$ and $t\in [0,1]$ such that $f((1-t)x+ty)<\min\{f(x),f(y)\}$. Let
\[T^*=\arg\min\{f((1-t)x+ty)|t\in [0,1]\}.\]
Then, $T^*$ is a closed subset of $[0,1]$ such that $0,1\notin T^*$. Let $t^*=\max T^*$ and $z=(1-t^*)x+t^*y$. Moreover, define $x(t)=(1-t)x+ty$. Because $Df(z)\neq 0$, we have that there exists $p\in X$ such that $\langle p,Df(z)\rangle>0$. Consider the following function $u$:
\[u(a,b)=f(z+a(y-x)+bp)=f(x(t^*+a)+bp).\]
Then, $u$ is $C^1$ around $(0,0)$, $u(0,0)=f(z)$, and $\frac{\partial u}{\partial b}(0,0)=\langle p,Df(z)\rangle>0$. By the implicit function theorem, there exist $\varepsilon_1, \varepsilon_2>0$, and a $C^1$ function $b:[-\varepsilon_1,\varepsilon_1]\to [-\varepsilon_2,\varepsilon_2]$ such that $b(0)=0$ and, for each $(a,b)\in [-\varepsilon_1,\varepsilon_1]\times [-\varepsilon_2,\varepsilon_2]$, $u(a,b)=f(z)$ if and only if $b=b(a)$. Then,
\[b'(a)=-\frac{\frac{\partial u}{\partial a}(a,b(a))}{\frac{\partial u}{\partial b}(a,b(a))}=-\frac{\langle y-x,Df(x(t^*+a)+b(a)p)\rangle}{\langle p,Df(x(t^*+a)+b(a)p)\rangle}.\]
Because of the definition of $z$ and the first-order condition, we have that
\[\langle y-x,Df(z)\rangle=0,\]
and thus,
\[b'(0)=-\frac{\langle y-x,Df(z)\rangle}{\langle p,Df(z)\rangle}=0.\]
Define
\[y(a)=x(t^*+a)+b(a)p,\ w(a)=b'(a)p+(y-x),\]
and
\[d(a)=\langle p,Df(y(a))\rangle.\]
Differentiating $u(a,b(a))=f(y(a))\equiv f(z)$ with respect to $a$, we obtain
\[\langle w(a),Df(y(a))\rangle=0,\]
and thus,
\[\langle w(a),g(y(a))\rangle=0\]
for any $a\in [-\varepsilon_1,\varepsilon_1]$. By our choice of $p$, $y(0)=z$ and $d(0)>0$, and thus, there exists $\delta>0$ such that if $0<a<\delta$, then $d(a)>0$. For such $a>0$,
\begin{align*}
0=&~\liminf_{a'\to a}\frac{1}{a'-a}[\langle w(a'),Df(y(a))\rangle-\langle w(a'),Df(y(a))\rangle]\\
=&~\liminf_{a'\to a}\frac{1}{a'-a}[\langle w(a'),Df(y(a))\rangle-\lambda(y(a))\langle w(a'),g(y(a))\rangle]\\
=&~\liminf_{a'\to a}\frac{1}{a'-a}[\langle w(a')-w(a),Df(y(a))\rangle+\lambda(y(a))\langle w(a'),g(y(a'))-g(y(a))\rangle] \\
=&~d(a)\times \liminf_{a'\to a}\frac{b'(a')-b'(a)}{a'-a}+\lambda(y(a))\langle w(a),Dg(y(a))w(a)\rangle\\
\le&~d(a)\times \liminf_{a'\to a}\frac{b'(a')-b'(a)}{a'-a}.
\end{align*}
Therefore,
\[\liminf_{a'\to a}\frac{b'(a')-b'(a)}{a'-a}\ge 0,\]
which implies that
\[\liminf_{a'\downarrow a}\frac{b'(a')-b'(a)}{a'-a}\ge 0,\ \liminf_{a'\uparrow a}\frac{b'(a')-b'(a)}{a'-a}\ge 0.\]
Fix $a\in ]0,\delta[$, and define $h(s)=b'(s)a-b'(a)s$. Then, $h(a)=h(0)=0$, and thus, there exists $s^*$ such that $0<s^*<a$ and $h(s)$ attains either the maximum or minimum on $[0,a]$ at $s=s^*$. If $h(s^*)$ attains the maximum, then
\begin{align*}
0\ge&~\liminf_{s\downarrow s^*}\frac{h(s)-h(s^*)}{s-s^*}\\
=&~a\liminf_{s\downarrow s^*}\frac{b'(s)-b'(s^*)}{s-s^*}-b'(a)\ge -b'(a),
\end{align*}
which implies that $b'(a)\ge 0$ for any $a\in ]0,\delta[$. If $h(s^*)$ attains the minimum, then
\begin{align*}
0\ge&~\liminf_{s\uparrow s^*}\frac{h(s)-h(s^*)}{s-s^*}\\
=&~a\liminf_{s\uparrow s^*}\frac{b'(s)-b'(s^*)}{s-s^*}-b'(a)\ge -b'(a),
\end{align*}
which again implies that $b'(a)\ge 0$ for any $a\in ]0,\delta[$. Because $b(0)=0$, we have that $b(a)\ge 0$ for all $a\in [0,\delta]$.

Because
\[\frac{\partial u}{\partial b}(0,0)=\langle p,Df(z)\rangle>0,\]
there exists a neighborhood $V$ of $(0,0)$ such that $\frac{\partial u}{\partial b}(a,b)>0$ for all $(a,b)\in V$. If $a>0$ is sufficiently small, then $a<\delta$ and $(a,b)\in V$ for all $b\in [0,b(a)]$. Therefore, $\frac{\partial u}{\partial b}(a,b)>0$ for such $(a,b)$, and thus,
\[f(x(t^*+a))\le f(x(t^*+a)+b(a)p)=f(z),\]
which contradicts the definition of $t^*$. This completes the proof. $\blacksquare$

\vspace{12pt}
Regarding strict quasi-concavity, we present the following theorem.

\vspace{12pt}
\noindent
{\bf Theorem 2}. Suppose that $U$ is an open and convex subset of a Banach space $X$ and $f:U\to \mathbb{R}$ is $C^1_*$. Choose a pair of a $C^1$ function $g$ and a positive continuous function $\lambda$ that satisfies (\ref{VEC}). If $\langle w,Dg(x)w\rangle<0$ for all $x\in U$ and $w\in X$ such that $w\neq 0$ and $\langle w,g(x)\rangle=0$, then $f$ is strictly quasi-concave.

\vspace{12pt}
\noindent
{\bf Proof}. Suppose that $f(x)$ is not strictly quasi-concave. Then, there exist $x,y\in U$ and $t\in ]0,1[$ such that $x\neq y$ and
\[f((1-t)x+ty)\le \min\{f(x),f(y)\}.\]
Define $x(t)=(1-t)x+ty$. Then, there exists $t^*\in \arg\min\{f(x(t))|t\in [0,1]\}$ such that $0<t^*<1$. Define $z=x(t^*)$, $w=y-x$, and
\[c(t)=\langle w,Df(x(t))\rangle.\]
By the first-order condition, we have that $c(t^*)=0$, and thus $\langle w,g(x(t^*))\rangle=0$. Moreover,
\begin{align*}
0>&~\langle w,Dg(z)w\rangle=\lim_{t\downarrow t^*}\frac{\langle w,g(x(t))-g(x(t^*))\rangle}{t-t^*}\\
=&~\lim_{t\downarrow t^*}\frac{\langle w,g(x(t))\rangle}{t-t^*}=\lim_{t\downarrow t^*}\frac{1}{\lambda(x(t))}\frac{\langle w,Df(x(t))\rangle}{t-t^*}\\
\ge&~\limsup_{t\downarrow t^*}\frac{1}{m}\frac{\langle w,Df(x(t))\rangle}{t-t^*}=\frac{1}{m}\limsup_{t\downarrow t^*}\frac{c(t)}{t-t^*},
\end{align*}
where $m=\min\{\lambda(x(t))|t\in [0,1]\}>0$. This implies that $c(t)<0$ for all $t>t^*$ such that $t-t^*$ is sufficiently small, and thus,
\[f(x(t))<f(z)\]
for such $t$, which is a contradiction. This completes the proof. $\blacksquare$

\section{Remarks}
Suppose that $X=\mathbb{R}^n$ and $U$ is an open and convex subset of $X$. Let $f:U\to \mathbb{R}$ be $C^1$. Then, $Df(x)$ can be represented by the row vector
\[\left(\frac{\partial f}{\partial x_1}(x),...,\frac{\partial f}{\partial x_n}(x)\right).\]
That is, we can consider that $Df(x)=(\nabla f(x))^T$. Moreover, suppose that $f$ is $C^2$. Then, $D^2f(x)$ can be represented by the following {\bf Hessian matrix}:
\[\begin{pmatrix}
\frac{\partial^2f}{\partial x_1^2}(x) & ... & \frac{\partial^2f}{\partial x_n\partial x_1}(x) \\
\vdots & \ddots & \vdots \\
\frac{\partial^2f}{\partial x_1\partial x_n}(x) & ... & \frac{\partial^2f}{\partial x_n^2}(x)
\end{pmatrix}.\]
In this context, it is known the following facts. Suppose that $f:U\to \mathbb{R}$ is $C^2$. Then, the following holds.

\begin{enumerate}[1)]
\item If $Df(x)\neq 0$ for all $x\in U$, then $f$ is quasi-concave if and only if $v^TD^2f(x)v\le 0$ for every $x\in U$ and $v\in \mathbb{R}^n$ such that $Df(x)v=0$.

\item If $Df(x)\neq 0$ and $v^TD^2f(x)v<0$ for every $x\in U$ and $v\in \mathbb{R}^n$ such that $v\neq 0$ and $Df(x)v=0$, then $f$ is strictly quasi-concave.
\end{enumerate}

Fact 1) is proved by Otani (1983). Fact 2) is well known, and have been introduced in many textbooks. If $Df(x)=0$ is admitted, then 1) is violated by the following function
\[f(x_1,x_2)=x_1^4.\]
Indeed, we have that
\[Df(x)=(4x_1^3,0), D^2f(x)=\begin{pmatrix}
12x_1^2 & 0\\
0 & 0
\end{pmatrix},\]
and thus $Df(x)v=0$ if and only if $x_1v_1=0$. For such a $v$, $v^TD^2f(x)v=0$, although $f$ is not quasi-concave.

Note that the converse of 2) is not true. Let us check this. Suppose that $U=\{(x_1,x_2)|x_1>0,\ x_2>0\}$. Define $f(x_1,x_2)=x_1^3x_2+x_1x_2^3$. If $x_2(x_1)$ represents a function whose graph consists of a level set $f(x_1,x_2)\equiv a$, then we can show that the function $x_2(x_1)$ is strictly convex, which implies that $f$ is strictly quasi-concave. We also have
\[Df(1,1)=(4,4),\ D^2f(1,1)=\begin{pmatrix}
6 & 6\\
6 & 6
\end{pmatrix},\]
and thus, if $v=(1,-1)$, then $Df(1,1)v=0$ but $v^TD^2f(1,1)v=0$. This example was found by Katzner (1968), and has an important implication in the consumer's optimization problem.

Actually, facts 1) and 2) are corollaries of our results. Recall that if $f$ is $C^2$ and $Df(x)\neq 0$ for all $x\in U$, then $f$ is $C^1_*$, and we can use $Df(x)$ as $g(x)$. Then, 1) follows from Theorem 1, and 2) follows from Theorem 2. In this view, our results are extensions of the above facts.

Is there an $f$ that is $C^1_*$ but is not $C^2$? we can answer this question affirmatively. Let
\begin{equation}
f(x)=\begin{cases}
x_2 & \mbox{if }x_2\le 0,\\
\frac{x_2}{1-x_1x_2} & \mbox{if }x_2>0.
\end{cases}\label{G0}
\end{equation}
Then, $f$ is $C^1$. Later, we will show that $f$ is not $C^2$. Define
\begin{align}
g_1(x)=&~\begin{cases}
0 & \mbox{if }x_2\le 0,\\
\frac{x_2^2}{\sqrt{1+x_2^4}} & \mbox{if }x_2>0,
\end{cases}\label{G1}\\
g_2(x)=&~\begin{cases}
1 & \mbox{if }x_2\le 0,\\
\frac{1}{\sqrt{1+x_2^4}} & \mbox{if }x_2>0,
\end{cases}\label{G2}\\
\lambda(x_1,x_2)=&~\begin{cases}
1 & \mbox{if }x_2\le 0,\\
\frac{\sqrt{1+x_2^4}}{(1-x_1x_2)^2} & \mbox{if }x_2>0.
\end{cases}\label{G3}
\end{align}
Then, we can easily check that $g$ is $C^1$, $\lambda$ is continuous, and $g,\lambda$ satisfy (\ref{VEC}).\footnote{The notation $Df(x)$ in (\ref{VEC}) is changed to $\nabla f(x)$ in this context.} Thus, $f$ is $C^1_*$. Note that, on some convex neighborhood of $0$, $f$ is quasi-concave; we can use Theorem 2 to check this fact.

It may be worthwhile to mention a fact on the above example. The function $g$ defined by (\ref{G1}) and (\ref{G2}) was first introduced by Debreu (1976).\footnote{In this paper, Debreu said that this example was first found by Andreu Mas-Colell and Leonard Shapiro, and Marcel Richter informed him of this example.} Debreu stated that, for this $g$, there is no $C^2$ function $f$ that satisfies (\ref{VEC}) with some positive continuous function $\lambda$ around $0$. Actually, this claim is correct. We now prove this fact rigorously.

First, we show that the function $f$ defined by (\ref{G0}) is not $C^2$ on any open neighborhood of $0$. If $f$ is $C^2$, then
\[\lambda(x)=\frac{\frac{\partial f}{\partial x_2}(x)}{g_2(x)},\]
which is $C^1$. However, if $x_1\neq 0$ and $x_2=0$, $\lambda$ defined by (\ref{G3}) is not differentiable, and thus $f$ is not $C^2$.

Next, suppose that there exists a $C^2$ function $h$ and a positive $C^1$ function $\mu$ defined on some open neighborhood of $0$ such that
\begin{equation}\label{VEC2}
\nabla h(x)=\mu(x)g(x)
\end{equation}
for all $x$. Without loss of generality, we can assume that the domain of $h$ includes $V=[-\varepsilon,\varepsilon]^2$, where $0<\varepsilon<1$. Then, $f$ is also defined on $V$. Consider the following differential equation:
\[\dot{x}_2(x_1)=-\frac{g_1(x_1,x_2(x_1))}{g_2(x_1,x_2(x_1))},\ x_2(0)=c,\]
where $-\varepsilon<c<\varepsilon$. By the Picard-Lindel\"of theorem, the nonextendable solution of the above equation is unique. Note that, by construction,
\[\dot{x}_2(x_1)=-\frac{\frac{\partial f}{\partial x_1}(x_1,x_2(x_1))}{\frac{\partial f}{\partial x_2}(x_1,x_2(x_1))}=-\frac{\frac{\partial h}{\partial x_1}(x_1,x_2(x_1))}{\frac{\partial h}{\partial x_2}(x_1,x_2(x_1))},\]
and thus, the graph of $x_2(x_1)$ coincides with the level sets $f^{-1}(c)$ and $h^{-1}(h(0,c))$.\footnote{Note that $f(0,c)=c$ for all $c$.} This implies that, on $V$, $h(x)=\varphi(f(x))$, where $\varphi(c)=h(0,c)$. By definition, we have that $\varphi$ is $C^2$. Moreover, by (\ref{VEC2}), we have that $\varphi'(c)\neq 0$, and
\[\frac{\partial h}{\partial x_2}(x_1,x_2)=\varphi'(f(x_1,x_2))\frac{\partial f}{\partial x_2}(x_1,x_2).\]
Because $h$ and $\varphi$ are $C^2$, we conclude that $\frac{\partial f}{\partial x_2}(x_1,x_2)$ is differentiable. However, we have already stated that it is not differentiable if $x_1\neq 0$ and $x_2=0$, which is a contradiction. Hence, such an $h$ does not exist.

There are two open problems related to our theorems. First, consider $X=\mathbb{R}^n$. Suppose that $A=(a_{ij})_{i,j=1}^n$ is an $n\times n$ symmetric matrix and $b\in \mathbb{R}^n$ such that $b_1\neq 0$. Debreu (1952) showed that $w^TAw\le 0$ for all $w\in \mathbb{R}^n$ such that $w^Tb=0$ if and only if for any permutation $\pi$ and $j\in \{2,...,n\}$,
\[(-1)^j\begin{vmatrix}
a_{i_1i_1} & ... & a_{i_1i_j} & b_{i_1}\\
\vdots & \ddots & \vdots & \vdots\\
a_{i_ji_1} & ... & a_{i_ji_j} & b_{i_j}\\
b_{i_1} & ... & b_{i_j} & 0
\end{vmatrix}\ge 0,\]
where $i_k=\pi(k)$. If the pair $(A,b)$ has this property, we say that $(A,b)$ satisfies {\bf Property N}. Applying this result to fact 1), we obtain the following result: suppose that $U\subset \mathbb{R}^n$ is an open and convex set, $f;U\to \mathbb{R}$ is $C^2$, and $\frac{\partial f}{\partial x_1}(x)\neq 0$ for all $x\in U$. Then, $f$ is quasi-concave if and only if, for each $x\in U$, the pair $(D^2f(x),\nabla f(x))$ satisfies Property N. We want to obtain the analogy of this fact. Suppose that $U\subset \mathbb{R}^n$ is open and convex, $f:U\to \mathbb{R}$ is $C^1_*$, and $g,\lambda$ satisfy (\ref{VEC}), where $g_1(x)\neq 0$ for all $x\in U$. Then, we expect that $f$ is quasi-concave if and only if, for each $x\in U$, the pair $(Dg(x),g(x))$ satisfies Property N. However, $Dg(x)$ may not be symmetric, and thus Debreu's result cannot be applied directly. Hence, whether our conjecture is true is an open problem.

Second, we would like to link our results to the {\bf concavity} of functions. Suppose that $U\subset X$ is open and convex, and $f:U\to \mathbb{R}$ is given. Then, $f$ is said to be {\bf concave} if and only if, for every $x,y\in U$ and $t\in [0,1]$,
\[f((1-t)x+ty)\ge (1-t)f(x)+tf(y).\]
If, in addition,
\[f((1-t)x+ty)>(1-t)f(x)+tf(y)\]
when $x\neq y$ and $0<t<1$, then $f$ is said to be {\bf strictly concave}. If $X=\mathbb{R}^n$ and $f$ is $C^2$, then it is well known that 1) $f$ is concave if and only if $D^2f(x)$ is negative semi-definite and 2) if $D^2f(x)$ is negative definite, then $f$ is strictly concave. As an analogy of these facts, we expect that for a $C^1_*$ function $f$, $f$ is concave if and only if $Dg(x)$ is negative semi-definite, and if $Dg(x)$ is negative definite, then $f$ is strictly concave. However, this conjecture is also difficult to verify. This is another open problem.

\section*{Acknowledgments}
The author is grateful to Shinichi Suda and Ken Hasegawa for their kindful comments.

\section*{Reference}

\begin{description}

\item{[1]} Aliprantis, C. D. and Border, K. C. (2007) Infinite Dimensional Analysis: A Hitchhiker's Guide. 3rd ed. Springer, Berlin.

\item{[2]} Debreu, G. (1952) ``Definite and Semi-Definite Quadratic Forms.'' Econometrica 20, pp.295-300.

\item{[3]} Debreu, G. (1972) ``Smooth Preferences.'' Econometrica 40, pp.603-615.

\item{[4]} Debreu, G. (1976) ``Smooth Preferences, A Corrigendum.'' Econometrica 44, pp.831-832.

\item{[5]} Katzner, D. W. (1968) ``A Note on the Differentiability of Consumer Demand Functions.'' Econometrica 36, pp.415-418.

\item{[6]} Luenberger, D. G. (1969) Optimization by Vector Space Methods. Wiley, New York.

\item{[7]} Otani, K. (1983) ``A Characterization of Quasi-Concave Functions.'' Journal of Economic Theory 31, pp.194-196.

\end{description}

\end{document}